# New applications of the Ahlfors Laplacian

Sergey E. Stepanov[1], Josef Mikeš[2], Irina I. Tsyganok[3]

[1] Department of Mathematics, Russian Institute for Scientific
and Technical Information of the Russian Academy of Sciences,
20, Usievicha street, 125190 Moscow, Russia,
E-mail address: s.e.stepanov@mail.ru

[2] Palacky University, 17. Listopadu 12, 77146 Olomouc,
Czech Republic
E-mail address: josef.mikes@upol.cz

[3] Department of Data Analysis and Financial Technologies,
Finance University, 49-55, Leningradsky Prospect,
125468 Moscow, Russia,
E-mail address: i.i.tsyganok@mail.ru

**Abstract.** In this article, we consider an orthogonal decomposition of the traceless part of the Ricci tensor of a compact Riemannian manifold and study its application to the geometry of compact almost Ricci solitons. In addition, we consider an orthogonal expansion of the traceless part of the second fundamental form of a compact spacelike hypersurface in a Lorentzian manifold and study its application to the problem of constructing solutions of general relativistic constraint equations in vacuum. In these two cases, we use the well-known Ahlfors Laplacian.

## 1   Introduction

The second-order elliptic operator $S^*S$ was considered over differential one-forms on an $n$-dimensional $(n \geq 2)$ compact Riemannian manifold $(M, g)$ without boundary in a number of papers [1] - [6]. The kernel of $S^*S$ is the space of conformal Killing one-forms. One of main applications of $S^*S$ is the construction of families of smooth quasi-conformal deformations of transformations $(M, g)$. The operator $S^*S$ was called as *the Ahlfors Laplacian* (see, for example, [4]; [6]). It is a particular form of the *Tachibana rough Laplacian* that is defined on differential $p$-forms for $1 \leq p \leq n-1$. In turn, the kernel of it is the space of conformal Killing $p$-forms. The geometry of the Tachibana rough Laplacian is described in [7]; [8] and etc.

On the other hand, one of the important components of the theory of Ricci flow are self-similar solutions called Ricci solitons (see [9, pp. 153-176] and [10, Chapter 3]). This theory, besides being known after G. Perelman's proof of the Poincar´e conjecture (for details see [10, Introduction]), has a wide range of applications in differential geometry and theoretical physics. In [11], was defined an almost Ricci soliton. This definition was formal and had no geometric meaning. In later works, only compact almost Ricci solitons were studied (e.g., [12], [13], [14] and etc.).

In the present paper we consider an application of the Ahlfors Laplacian to the geometry of almost Ricci solitons. First, this allows us to substantiate the concept of almost Ricci solitons. Second, the statements proved earlier become corollaries of the corresponding properties of this orthogonal decomposition. The purpose of the fourth section is to establish some algebraic-differential properties of the general relativistic vacuum constraint equations due to the associated Ahlfors Laplacian.

## 2. The orthogonal decomposition of the traceless Ricci tensor on a compact Riemannian manifold

Let $(M, g)$ be an $n$-dimensional $(n \geq 3)$ compact Riemannian manifold with the Levi-Civita connection $\nabla$. We denote by $S^p M := S^p T^* M$ the vector space of symmetric bilinear differential $p$-forms on $M$, or, in other words, of covariant symmetric $p$-tensors $(p \geq 1)$ on $(M, g)$ and define the $L^2$ *inner product* of two covariant symmetric $p$-tensors $\varphi$ and $\varphi'$ on $(M, g)$ by the formula

$$\langle \varphi, \varphi' \rangle := \int_M g(\varphi, \varphi') \, dvol_g$$

where $dvol_g$ being the volume element of $(M, g)$. Also $\delta^*: C^\infty S^1 M \to C^\infty S^2 M$ will be the first-order differential operator defined by the formula $\delta^* \theta := \frac{1}{2} L_\xi g$ for some smooth vector field $\xi$ and it's $g$-dual one-form $\theta$ (see [15, p. 117; 514]). At the same time, we denote by the formula $\delta: C^\infty S^2 M \to C^\infty S^1 M$ the formal adjoint operator for $\delta^*$ which is called the *divergence of symmetric two-tensors*. In this case, we have $\langle \varphi, \delta^* \theta \rangle = \langle \delta \varphi, \theta \rangle$ for any $\varphi \in C^\infty S^2 M$ and $\theta \in C^\infty S^1 M$. Using the above, we can definite the $S$ from $S^1 M$ to the vector space of trace free symmetric tensors $S_0^2 M$ defined by the equality $S\theta := \delta^*\theta + \frac{1}{n}\delta\theta\, g$. The operator $S$ is known as the *Cauchy-Ahlfors* operator (see, for example, [3]). It's obvious that a vector field $\xi$ on $(M,g)$ is a *conformal Killing vector field* if and only if $S$ annihilates the one-form $\theta$ identified with $\xi$ by the Riemannian metric $g$, i.e., $\theta^{\#} = \xi$ (see [16, p. 46]). Therefore, the kernel of $S$ is the finite-dimensional space.

The formal adjoint operator for $S$ is defined by the formula $S^*\omega = \delta\omega$ for an arbitrary $\omega \in C^\infty S_0^2 M$ (see [4]; [5]). Moreover, it is easy to verify the formula (see also [3]; [4] and [5])

$$\triangle_A \theta = \frac{1}{2}\delta d\theta - Ric(\xi,\cdot) + \frac{n-1}{n} d\delta\theta. \tag{2.1}$$

for the *Ahlfors Laplacian* $\triangle_A := S^*S$. In (2.1), we denote by $Ric$ the Ricci tensor of $(M,g)$ and by $d$ the well known exterior derivative from the vector space of differential one-forms $\Lambda^1 M$ to the vector space of differential two-forms $\Lambda^2 M$ (see [8, p. 21]). We recall that $\triangle_A$ is a formally self-adjoint non-negative strongly elliptic second order differential operator (see [2]; [3]; [5]). Note that $ker\,\triangle_A = ker\,S$ since $\langle\theta, \triangle_A \theta\rangle = \langle\theta, S^*S\theta\rangle = \langle S\theta, S\theta\rangle$ for any $\theta \in C^\infty S^1 M$. Hence the kernel of $\triangle_A$ consists of conformal Killing one-forms.

We recall here that a symmetric divergence free and traceless covariant two-tensor is called $TT$-tensor (see, for instance, [17]). Any $TT$-tensor is denoted by $\varphi^{TT}$ (see [19]). As a consequence of a result of Bourguignon-Ebin-Marsden (see [15, p. 132]; [19]) the space of $TT$-tensors is an infinite-dimensional vector space on any compact Riemannian manifold $(M,g)$. Such tensors are of fundamental importance in stability analysis in General Relativity (see, for instance, [18]) and in Riemannian geometry (see, for instance, [15, p. 346-347]; [19] and [20]).

Now, we are ready to formulate our first result.

**Theorem 1.** *Let $(M,g)$ be an $n$-dimensional $(n \geq 3)$ compact Riemannian manifold. Then there are some one-form $\theta \in C^\infty S^1 M$ and some $TT$-tensor $\varphi^{TT} \in C^\infty S^2 M$ such that the traceless part $Ric_0$ of the Ricci tensor $Ric$ of $(M,g)$ has the decomposition*

$$Ric_0 = S\theta + \varphi^{TT} \tag{2.2}$$

*for the Cauchy-Ahlfors operator S. Furthermore, both factors in the right side of (2.2) are orthogonal to each other with respect to the $L^2$ inner product and the following equation holds*

$$\triangle_A \theta = -\frac{n-2}{n} ds \qquad (2.3)$$

*for the Ahlfors Laplacian $\triangle_A$.*

*Proof.* For any $n$-dimensional ($n \geq 3$) compact Riemannian manifold $(M, g)$, the algebraic sum $\text{Im } \delta^* + C^\infty M \cdot g$ is closed in $S^2 M$, and we have the decomposition

$$S^2 M = (\text{Im } \delta^* + C^\infty M \cdot g) \oplus \left(\delta^{-1}(0) \cap \text{trace}_g^{-1}(0)\right) \qquad (2.3)$$

where both factors are infinite dimensional and orthogonal to each other with respect to the $L^2$ inner scalar product (see Lemma 4.57 in [15, p. 130]). It's obvious that the second factor $\delta^{-1}(0) \cap \text{trace}_g^{-1}(0)$ of (2.4) is the space of *TT*-tensors. In particular, it is well known that $Ric \in S^2 M$ for the Ricci tensor $Ric$ of $(M, g)$. Therefore, in the case of the Ricci tensor, from orthogonal decomposition (2.4) we obtain

$$Ric = \left(\frac{1}{2} L_\xi g + \lambda g\right) + \varphi^{TT} \qquad (2.5)$$

for some vector field $\xi \in C^\infty TM$, some *TT*-tensor $\varphi^{TT} \in C^\infty S^2 M$ and some scalar function $\lambda \in C^\infty M$. It is obvious that the vector field $\xi$ is defined up to a *Killing vector field*. We recall here that any such vector field $X$ satisfies the equation $L_X g = 0$ (see [16, p. 43]).

We denote by $s := \text{trace}_g Ric$ the scalar curvature of $(M, g)$. Then applying the operator $\text{trace}_g$ to both sides of (2.5), we obtain

$$s = -\delta\theta + n\lambda \qquad (2.6)$$

where $\theta$ is the $g$-dual one-form of $\xi$. In this case, (2.4) can be rewritten in the form (2.2) where $Ric_0 = Ric - \frac{1}{n} s\, g$ is the traceless part of the Ricci tensor $Ric$ and $S\theta = L_\xi g + \frac{1}{n} \delta\theta\, g$ denotes the Cauchy-Ahlfors operator actions on the one-form $\theta$. Next, applying $S^*$ to both sides of (2.2), we obtain (2.3), because of $\delta\, Ric = -\frac{1}{2} ds$, which is a consequence of the second Bianchi identity (see [15, p. 43]). The proof of the theorem is complete.

The following corollary holds.

**Corollary 1.** *Let $(M, g)$ be an $n$-dimensional $(n \geq 3)$ compact Riemannian manifold and $Ric_0 = S\theta + \varphi^{TT}$ be a $L^2$ orthogonal decomposition of the traceless Ricci tensor $Ric_0$ of $(M, g)$ for some one-form $\theta \in C^\infty S^1 M$ and some TT-tensor $\varphi^{TT} \in C^\infty S^2 M$. If one of the following assumptions $\int_M L_\xi s \, dv_g = 0$ for $\xi = \theta^\#$ or $s = constant$ holds, then $\xi$ is a conformal Killing vector field and $\varphi^{TT}$ is the traceless Ricci tensor $Ric_0$.*

*Proof.* From (2.3) we deduce the integral formula

$$\langle S\theta, S\theta \rangle = -\frac{n-2}{n} \int_M L_\xi s \, dv_g \qquad (2.7)$$

Now, if one of the two conditions $\int_M L_\xi s \, dv_g = 0$ or $s = constant$ holds, then from (2.7) we obtain $S\theta = 0$, i.e., $\xi = \theta^\#$ is a conformal Killing vector field. Then (2.2) can be rewritten in the form $Ric - \frac{s}{n} g = \varphi^{TT}$ where $s = constant$. In this case, $\varphi^{TT}$ is the traceless part of the Ricci tensor $Ric_0 = Ric - \frac{s}{n} g$ of $(M, g)$ for $s = constant$.

**Remark.** In particular, from Corollary 1 we conclude that the scalar curvature $s$ of a compact Riemannian manifold $(M, g)$ is constant if and only if the Ricci tensor $Ric$ of $(M, g)$ has the form $Ric = \frac{s}{n} g + \varphi^{TT}$ for some TT-tensor $\varphi^{TT} \in C^\infty S^2 M$.

We recall that if a Riemannian $n$-dimensional manifold $(M, g)$ admits a complete conformal Killing vector field $\xi$ with singular points at each of which its divergence $div \, \xi$ does not vanish, then $(M, g)$ is homeomorphic to either the standard sphere $\mathbb{S}^n$ or the Euclidian space $\mathbb{R}^n$ (see Lemma 2.2 in [21]). On the other hand, if $(M, g)$ is a compact manifold, then vector field $\xi$ over $(M, g)$ is complete. Then using Lemma 2.2 from [21] and our Corollary 1, we can formulate the following proposition.

**Corollary 2.** *Let $(M, g)$ be an $n$-dimensional $(n \geq 3)$ compact Riemannian manifold and $Ric_0 = S\theta + \varphi^{TT}$ be a $L^2$ orthogonal decomposition of its traceless Ricci tensor $Ric_0$ for some TT-tensor $\varphi^{TT} \in C^\infty S^2 M$ and some one-form $\theta \in C^\infty S^1 M$ such that $\xi = \theta^\#$ has singular points at each of which $\delta \theta \neq 0$. Then from*

*one of the assumptions $\int_M L_V s \, dv_g = 0$ or, in particular, $s = constant$ one can deduce that $(M, g)$ is locally homeomorphic the standard sphere $\mathbb{S}^n$ or the Euclidian space $\mathbb{R}^n$.*

We finish with a remark. Namely, we note that two special classes of Riemannian manifolds can be defined based on the orthogonal decomposition (2.2). The first class is defined by the condition $S\theta = 0$. In this case, any such Riemannian manifold satisfies the equality $Ric = \frac{s}{n} g + \varphi^{TT}$ for same $TT$-tensor $\varphi^{TT}$. In turn, the second class is defined by the condition $\varphi^{TT} = 0$. In this case, any such Riemannian manifold satisfies the equation $Ric = \frac{1}{2} L_\xi g + \lambda g$. Any such Riemannian manifold is called an almost Ricci soliton. We will consider this class of Riemannian manifolds in the next section.

## 3. Compact almost Ricci solitons

An $n$-dimensional $(n \geq 2)$ Riemannian manifold $(M, g)$ is said to be a Ricci (respectively, an *almost Ricci*) soliton (see [9, pp. 37-38] and [11], respectively) if there exist a smooth vector field $V$ and a constant (respectively, a smooth function) $\lambda$ such that

$$Ric = \frac{1}{2} L_V g + \lambda g. \quad (3.1)$$

It is obvious that the vector field $V$ is defined up to a Killing vector field. We denote by $(M, g, V, \lambda)$ the Ricci (respectively, an *almost Ricci*) soliton, and (3.1) we call its structural equations. We notice here that when $n \geq 3$ and $V$ is a *Killing vector field* an almost Ricci soliton $(M, g, V, \lambda)$ will be a Ricci soliton, since in this case we have an *Einstein manifold*, that is, satisfying $Ric = \frac{s}{n} g$ where $s$ is constant. If, in addition, $\xi \equiv 0$, then an almost Ricci soliton $(M, g, V, \lambda)$ is called *trivial*.

Applying the operator $trace_g$ to both sides of (3.1), we obtain $s = -\delta\theta + n\lambda$. In this case we can rewrite Equation (3.1) in the form $Ric_0 = S\theta$ for $V = -\theta^\#$. Therefore, the following corollary holds.

**Corollary 3.** *Let $(M, g, V, \lambda)$ be an $n$-dimensional $(n \geq 3)$ compact almost Ricci soliton with non-constant $\lambda$. If its scalar curvature $s$ satisfies one of the assumptions*

$\int_M L_V s \, dv_g = 0$ *or, in particular,* $s = constant$, *then* $(M, g, V, \lambda)$ *is isometric to the standard sphere* $\mathbb{S}^n$ *and* $V$ *is a gradient conformal Killing vector field.*

*Proof.* Let $(M, g, V, \lambda)$ be an $n$-dimensional ($n \geq 3$) compact almost Ricci soliton. If one of the conditions $\int_M (L_V s) \, dvol_g = 0$ or $s = constant$ is satisfied, then from Corollary 1 we conclude that $(M, g, V, \lambda)$ is an Einstein manifold. In this case, (3.1) can be rewritten in th form $\frac{1}{2} L_V g = \left(\lambda - \frac{s}{n}\right) g$. If $\lambda$ is a non-constant function, then $V$ is a non-isometric conformal Killing vector field on a compact Einstein manifold $(M, g, V, \lambda)$. In this case, $(M, g, V, \lambda)$ is isometric to the standard sphere $\mathbb{S}^n$ (see [16, p. 57]). On the other hand, if a sphere $\mathbb{S}^n$ of dimension $n \geq 2$ admits a globally defined conformal Killing vector field $Z$, then $Z$ is globally decomposed into the form $Z = X + Y$, where $X$ is a Killing vector field and $Y$ is a gradient conformal Killing vector field (see Theorem 3.1 in [24]). In turn, $V$ is a conformal Killing vector field of $(M, g, V, \lambda)$, then it is defined up to a Killing vector field. Therefore, we can say that $V$ is a gradient conformal Killing vector field. In turn, a sphere $\mathbb{S}^n$ of dimension $n \geq 2$ admits a globally defined gradient conformal Killing vector field (see Remark 3.2 in [24]). That finishes the proof of Corollary 3.

It has been known for some time that compact two-dimensional Ricci solitons must have constant curvature (see Theorem 10.1 in [21]). In turn, we can formulate the following statement on compact almost Ricci soliton for the case when $\dim M = 2$.

**Corollary 4.** *Let* $(M, g, V, \lambda)$ *be a two-dimensional compact almost Ricci soliton with nonpositive scalar curvature* $s$. *Then*

(i) *if* $s < 0$, *it is a trivial soliton;*

(ii) *if* $V$ *is a nowhere vanishing vector field, it is locally isometric to the Euclidean plane* $\mathbb{R}^2$.

*Proof.* Let $(M, g, V, \lambda)$ be a two-dimensional compact almost Ricci soliton, then $Ric = \frac{s}{2} g$ (see [15, p. 44]). First, from (3.1) we obtain the differential equation

$$ds = -d\delta\theta + n \, d\lambda. \tag{3.2}$$

Second, applying the operator $\delta$ to both sides of (3.1), we get the following equation:

$$ds = -\overline{\triangle}\,\theta - d\delta\theta + \frac{s}{2}g + n\,\lambda. \tag{3.3}$$

Combining (3.2) and (3.3) we deduce the second order differential equation

$$\overline{\triangle}\,\theta = \frac{s}{2}g. \tag{3.4}$$

In turn, from (2.5) we deduce the integral formula

$$\int_M s\,g(\xi,\xi)\,dv_g = 2\,\langle \nabla\xi, \nabla\xi \rangle \tag{3.5}$$

for $\xi = \theta^{\#}$. In this case, if the scalar curvature $s \leq 0$, then $\nabla\xi = 0$ and $\frac{s}{n}\|\xi\|^2 = 0$. Therefore, if $\xi$ is a nowhere vanishing vector field, then from the last equation we conclude that $s = 0$. Then $(M,g)$ has zero sectional curvature and hence $(M,g,V,\lambda)$ is locally isometric to the Euclidean plane $\mathbb{R}^2$ (see [7, p. 44]). On the other hand, if the scalar curvature $s \neq 0$ at any point $x \in M$, then $\xi = 0$ and hence $(M,g,V,\lambda)$ is trivial. We have proved our corollary.

## 4. A compact hypersurface in a Lorentzian manifold

Let $\mathbb{L}_1^{n+1}$ be a Lorentzian manifold of dimension $n+1 \geq 4$ and signature $(- + \cdots +)$ solving the Einstein field equations (see [26]; [27])

$$\overline{Ric} - \frac{1}{2}\bar{s}\,\bar{g} + 2\Lambda\,g = \kappa\,T \tag{3.1}$$

where we denote by $\bar{g}$ the metric tenor of $\mathbb{L}_1^{n+1}$, we also denote by $\overline{Ric}$ and $\bar{s}$ the Ricci curvature and the scalar curvature of $\bar{g}$, respectively. We employ the letter $\kappa$ for a positive constant, whose value (and physical dimensions) depends on the specific conventions one adopts. In addition, as it is rather customary in the physical literature, $T$ stands for the *stress-energy tensor* of the sources while $\Lambda$ stands for a *cosmological constant*. Any such $\mathbb{L}_1^{n+1}$ is called a *spacetime*.

Moreover, we consider a compact spacelike hypersurface $(M,g)$ in $\mathbb{L}_1^{n+1}$ (see, for example, [15, p. 38-39]). The extrinsic geometry of $(M,g)$ is described by the so-called *second fundamental form* $K \in S^2 M$. The second fundamental form $K$ has the following $L^2$ orthogonal decomposition

$$K = \left(\frac{1}{2}L_\xi g + \lambda\,g\right) + \varphi^{TT} \tag{3.2}$$

for some vector field $\xi \in C^\infty TM$, some TT-tensor $\varphi^{TT} \in C^\infty S^2 M$ and some scalar function $\lambda \in C^\infty M$.

We denote by $H := trace_g K$ the *mean curvature* of $(M, g)$ (see [15, p. 38]). Then applying the operator $trace_g$ to both sides of (3.2), we obtain

$$H = -\delta\theta + n\lambda \tag{3.3}$$

where $\theta$ is the $g$-dual one-form of $\xi$. Let $(M,g)$ be a compact manifold without boundary. In this case, (3.2) can be rewritten in the form

$$K_0 = S\theta + \varphi^{TT} \tag{3.4}$$

where $K_0 = K - \frac{1}{n} H\, g$ is the traceless part of the second fundamental form $K$ and $S\theta = L_\xi g + \frac{1}{n}\delta\theta\, g$ denotes the Cauchy-Ahlfors operator. Next, applying $S^*$ to both sides of (3.4), we obtain

$$\triangle_A \theta = \delta K + \frac{1}{n} dH \tag{3.5}$$

for the Ahlfors Laplacian $\triangle_A$. Now, we are ready to formulate our result.

**Theorem 2.** *Let $(M,g)$ be an $n$-dimensional compact spacelike hypersurface in a Lorentzian manifold $\mathbb{L}_1^{n+1}$ of dimension $n+1 \geq 4$. Then there are some one-form $\theta \in C^\infty S^1 M$ and some TT-tensor $\varphi^{TT} \in C^\infty S^2 M$ such that the traceless part $K_0$ of the second fundamental form $K$ of $(M,g)$ has the decomposition*

$$K_0 = S\theta + \varphi^{TT} \tag{3.6}$$

*for the Cauchy-Ahlfors operator $S$. Furthermore, both factors in the right side of (3.6) are orthogonal to each other with respect to the $L^2$ inner product and the following equation holds*

$$\triangle_A \theta = \delta K + \frac{1}{n} dH \tag{3.7}$$

*for the Ahlfors Laplacian $\triangle_A$.*

It is then a well-known fact that the triple $(M, g, K)$ solves a system of general relativistic constraint equations that takes the form (see [26]; [27])

$$\begin{cases} s - g(K,K) + (trace_g K)^2 = 2\, T(N,N) + 2\Lambda; \\ div_g K - d(trace_g K) = \kappa\, T(N, \cdot). \end{cases} \tag{3.8}$$

for a timelike unit normal vector field $N$ to $(M, g)$.

To avoid unnecessary complications, we shall focus here on the vacuum case, i.e., we consider the field equations with no sources ($T = 0$) and take $\Lambda = 0$. In this case, (3.1) can be rewritten in the form $\overline{Ric} = \frac{1}{2}\bar{s}\bar{g}$, then the field equations are also referred to as the vacuum field equations. By setting $T = 0$ and $\Lambda = 0$ in the trace-reversed field equations, the vacuum equations can be written as $Ric = 0$. Then $\mathbb{L}_1^{n+1}$ is a *Ricci-flat spacetime*.

**Remark.** The first example of this is the *Schwarzschild geometry* that describes a static black hole. In this geometry the Ricci tensor is zero everywhere but the spacetime is most certainly not flat. The second example is *Calabi-Yau manifolds* are Ricci flat, but they are certainly not the same as Minkowski space $\mathbb{R}_1^{n+1}$.

The well-known problem is to construct solutions of the general relativistic vacuum constraint equations (see [27]; [28]; [29] and etc.):

$$\begin{cases} s - g(K,K) + (trace_g K)^2 = 0; \\ div_g K - d(trace_g K) = 0. \end{cases} \quad (3.9)$$

These are the general relativistic constraint equations whatever the space-dimension $n \geq 2$. Moreover, we known the following existence theorem.

**Theorem 3** (see [27]). *Let $(M, g, K)$ be a triple satisfying the vacuum constraint equations (3.9), then there exists a spacetime $\mathbb{L}_1^{n+1}$ with metric $\bar{g}$ and an embedding $f: M \to \mathbb{L}_1^{n+1}$ such that the following assertions are true*:

1. *the spacetime $\mathbb{L}_1^{n+1}$ is Ricci-flat*;

2. *g is the induced metric by $f$, namely $g = f^*(\bar{g})$*;

3. *K is the second fundamental form of $f: M \to \mathbb{L}_1^{n+1}$*.

In turn, the second equation of (3.9) can be rewritten in the form

$$\delta K = -dH. \quad (3.10)$$

Combining (3.7) and (3.10), we obtain

$$\triangle_A \theta = -\frac{n-1}{n} dH.$$

Then the following corollary holds.

**Corollary 5.** *Let $(M, g)$ be an $n$-dimensional $(n \geq 3)$ compact spacelike hypersurface in a spacetime $\mathbb{L}_1^{n+1}$ and $K_0 = S\theta + \varphi^{TT}$ be a $L^2$ orthogonal*

*decomposition of the traceless part $K_0$ of its second fundamental form for some one-form $\theta \in C^\infty S^1 M$ and TT-tensor $\varphi^{TT} \in C^\infty S^2 M$. If the general relativistic vacuum constraint equations are true, then the mean curvature H of $(M, g)$ is constant if and only if $K = \frac{1}{n} H g + \varphi^{TT}$.*

We note that if the mean curvature $H$ of $(M, g)$ is constant and the second factor $\varphi^{TT}$ in (3.6) is zero, then $(M, g)$ is *totally umbilical* (see [15, p. 38]). On the other hand, we consider Lorentzian manifolds which admit a timelike conformal vector field. Any such manifold is called in [30] a *conformally stationary spacetime*. In turn, if $\mathbb{L}_1^{n+1}$ is a conformally stationary and Einstein spacetime, then its every compact spacelike hypersurface $(M, g)$ of constant mean curvature $H$ is totally umbilical (see also [30]). Therefore, we can formulate the following corollary.

**Corollary 6.** *Let $\mathbb{L}_1^{n+1}$ be a conformally stationary and Einstein spacetime for dimension $n + 1 \geq 4$ and let the general relativistic vacuum constraint equations are true. If $K = \left(S\theta + \frac{1}{n} Hg\right) + \varphi^{TT}$ is the $L^2$ orthogonal decomposition of the second fundamental form of a compact spacelike hypersurface $(M, g)$ with constant mean curvature, then its second factor $\varphi^{TT}$ is zero $(M, g)$.*

## References

1. Ahlfors L., Conditions for quasiconformal deformation in several variables, Contributions to Analysis. A Collection of Papers Dedicated to L. Bers, Academic Press, New York, 19-25 (1974).

2. Pierzchalski A., Ricci curvature and quasiconformal deformation of a Riemannian manifold, Manuscripta Math., **66**, 113-127 (1989).

3. Pierzchalski A., Gradients: the ellipticity and the elliptic boundary conditions – a jigsaw puzzle, Folia Mathematica, **19**(1), 65-83 (2017).

4. Pierzchalski A., Orsted B., The Ahlfors Laplacian on a Riemannian manifold with boundary, Michigan Math. J., **43**, 99-122 (1996).

5. Kozlowski W., Pierzchalski A., Natural boundary value problems for weighted form Laplacians, Ann. Scuola Norm. Sup. Pisa Cl. Sci., **VII** (5), 343-367 (2008).

6. Branson T., Stein-Weiss operators and ellipticity, Journal of Functional Analysis, **151**, 334-383 (1997).

7. Stepanov, S. E.; Mikesh, J. The Hodge–de Rham Laplacian and Tachibana operator on a compact Riemannian manifold with curvature operator of definite sign. (Russian) Izv. Ross. Akad. Nauk Ser. Mat., **79** (2), 167–180 (2015); translation in Izv. Math., **79** (2), 375–387 (2015).

8. Stepanov S.E., Mikeš J., Betti and Tachibana numbers of compact Riemannian manifolds, Differential Geometry and its Applications, **31**, 486-495 (2013).

9. Chow B., Lu P., and Ni L., Hamilton's Ricci flow, Grad. Stud. in Math., **77**, Amer. Math. Soc., Providence, RI (2006).

10. Morgan J., Tian G., Ricci flow and Poincare conjecture, Clay Mathematics Monographs, American Mathematical Society, Providence, RI; Clay Mathematics Institute, Cambridge, MA (2007).

11. Pigola S., Rigoli M., Rimoldi M., Setti A. G., Ricci almost solitons, Ann. Scuola Norm. Sup. Pisa Cl. Sci. **10** (4), 757–799 (2011).

12. Barros A., Batista R., Ribeiro jr. E., Compact almost Ricci solitons with constant scalar curvature are gradient, Monatshefte für Mathematik, **174**, 29–39 (2014).

13. Barros A., Gomes J.N., Rebeiro E., A note on rigidity of almost Ricci soliton, Archiv der Mathematik, **100**, 481–490 (2013).

14. Deshmukh S., Almost Ricci solitons isometric to spheres. Int. J. of Geom. Methods in Modern Physics, **16** (5), 1950073 (2019).

15. Besse A.L., Einstein manifolds, Springer-Verlag, Berlin & Heidelberg (2008).

16. Yano K., Integral formulas in Riemannian geometry, Marcel Dekker, New York (1970).

17. Gicquaud R., Ngo Q.A., A new point of view on the solutions to the Einstein constraint equations with arbitrary mean curvature and small *TT*-tensor, Class. Quant. Grav., **31**(19), 195014 (2014).

18. Garattini R., Self sustained tranversable wormholes?, Class. Quant. Grav., **22** (6), 2673–2682 (2005).

19. Bourguignon, J. P. and Ebin, David G. and Marsden, Jerrold E. Sur le noyau des opérateurs pseudo-differentiels á symbole surjectif et non injectif, Comptes rendus hebdomadaires des séances de l'Académie des sciences. Séries A et B, Sciences mathématiques et Sciences physiques, 282 (1976), 867-870.

20. Boucetta M., Spectra and symmetric eigentensors of the Lichnerowicz Laplacian on $S^n$, Osaka J. Math. **46** (2009) 235–254.

21. Obata M., Conformal transformations of Riemannian manifolds, J. Differential Geometry, **4** (3), 311-333 (1970).

22. Sharma R., Some results on almost Ricci solitons and geodesic vector fields, Beitrage zur Algebra und Geometrie, **59**, 289-294 (2018).

23. Ghosh, A., Ricci almost solitons satisfying certain conditions on the potential vector field, Publ. Math. Debrecen, **87** (1-2), 103-110 (2015).

24. Yano K., Sawaki S., Riemannian manifolds admitting a conformal transformation group, J. Differential Geometry, **2** (1968), 161-184.

25. Hamilton R.S., The Ricci flow on surfaces, Mathematics and general relativity, Contemp. Math., **71**, Amer. Math. Soc., Providence, 237–262 (1988).

26. Shoke-Bryua Y., Mathematical problems in general relativity, Uspekhi Mat. Nauk, **40** (6), 3–39 (1985).

27. Carlotto A., The general relativistic constraint equations, Living Reviews in Relativity, **24** (2), 1-170 (2021).

28. Chru P. T., Delay E., On mapping properties of the general relativistic constraints operator in weighted function spaces, with applications, Mémoires de la Société Mathématique de France, **94** (2003).

29. Helmut F., Cauchy problems for the conformal vacuum field equations in General Relativity, Commun. Math. Phys., **91**, 445-472 (1983).
30. Alías L.J., Romero A., Sanchez M., Spacelike hypersurfaces of constant mean curvature and Calabi-Bernstein type problems, Tôhoku Math J., **49**, 337-345, (1997).